\documentclass{amsart}

\usepackage[T1]{fontenc}
\usepackage{times}
\usepackage{amssymb,epsfig,verbatim}
\usepackage{calligra,mathrsfs}
\theoremstyle{plain}
\usepackage[T1]{fontenc}
\usepackage[utf8]{inputenc}
\usepackage[english]{babel}
\babeltags{en = english}
\babeltags{fr = frenchb}

\usepackage{amssymb}
\usepackage{amsthm}
\usepackage{graphics}
\usepackage{amsmath}
\usepackage{amstext}
\usepackage{multirow}
\usepackage{enumerate}

\usepackage[all,cmtip]{xy}

\usepackage{breqn}

\usepackage{hyperref}
\usepackage{graphicx}
\usepackage{amsmath}
\usepackage{amssymb}
\usepackage{amsthm}
\usepackage{amstext}
\usepackage{epstopdf}
\usepackage{mathrsfs}
\usepackage{amscd}
\usepackage{tikz}
\usetikzlibrary{cd}
\AtBeginDocument{%
   \def\MR#1{}
}

\providecommand{\keywords}[1]
{
  \small	
  \textbf{\textit{Keywords---}} #1
}

\newtheorem{thm}{Theorem}[section]

\newtheorem{theorema}{Theorem}

\newtheorem*{akn}{Acknowledgements}

\newtheorem{propositione}[theorema]{Proposition}

\newtheorem{dfn}[thm]{Definition}
\newtheorem{lemma}[thm]{Lemma}
\newtheorem{prop}[thm]{Proposition}

\newtheorem{cor}[thm]{Corollary}

\newtheorem{rmk}[thm]{Remark}

\newtheorem{ex}[thm]{Example}

\newcommand{\mb}{\mathbb}

\newcommand{\mc}{\mathcal}

\newcommand{\C}{\mb C}
\newcommand{\Pj}{\mb P}
\newcommand{\Z}{\mb Z}

\newcommand{\F}{\mc F}
\newcommand{\D}{\mc D}
\newcommand{\G}{\mc G}

\newcommand{\Oc}{\mathcal O}
\newcommand{\X}{X}
\newcommand{\Y}{Y}

\newcommand{\Qu}{Q}
\newcommand{\Fol}{\F\text{ol}}

\newcommand\restr[2]{{
  \left.\kern-\nulldelimiterspace 
  #1 
  \vphantom{\big|} 
  \right|_{#2} 
  }}

\DeclareMathOperator{\codim}{codim}
\DeclareMathOperator{\sing}{sing}

\DeclareMathOperator{\Aut}{Aut}
\DeclareMathOperator{\Pic}{Pic}

\DeclareMathOperator{\rest}{rest}

\DeclareMathOperator{\PGL}{PGL}
\DeclareMathOperator{\PO}{PO}

\DeclareMathOperator{\Hp}{\mathcal{H}}
\DeclareMathOperator{\Sp}{\mathcal{S}}

\DeclareMathOperator{\HomSheaf}{\mathscr{H}\text{\kern -3pt {\calligra\large om}}\,}

\numberwithin{equation}{section}

\addtocounter{section}{0}             
\numberwithin{equation}{section}       

\makeatletter
\def\keywords{\xdef\@thefnmark{}\@footnotetext}
\makeatother

\title{Extensions and restrictions of holomorphic foliations}
\author{Mateus Gomes Figueira}
\address{Instituto de Matemática Pura e Aplicada\\
Estrada Dona Castorina 110, 22460-320 Rio de Janeiro RJ, Brazil.}
\curraddr{IRMAR, Université de Rennes 1, Campus de Beaulieu, 35042 Rennes Cedex, France}
\email{mateus.gomes.figueira@hotmail.com}

\date{\today}

\begin{document}
\thanks{The author appreciates the financial support given by CNPq and Program CAPES/COFECUB.}

\begin{abstract}
 We prove an extension criterion for codimension one foliations on projective hypersurfaces based on the degree of the foliation and the degree of the hypersurface, and we ensure, in some instances, an isomorphism between the corresponding spaces of foliations. We also present some examples of foliations that do not satisfy the extension criterion and do not extend.
\end{abstract}

\maketitle
\setcounter{tocdepth}{1}
\tableofcontents

\section{Introduction}

Let $\F$ be a codimension one singular foliation on a smooth hypersurface $\X$ of $\Pj^n$, $n>2$. We say that a foliation $\G$ on $\Pj^n$ is an {\it extension} to $\F$ if its restriction to $\X$ is $\F$. D. Cerveau in \cite{MR3167127} proposed the investigation of necessary conditions to guarantee the existence of extensions of foliations on projective hypersurfaces. He further asks whether the unconditional existence of extensions for foliations on a hypersurface characterizes hyperplanes. In this work, we will show the following criterion for the existence of extensions of foliations on smooth projective hypersurfaces.

\begin{theorema}\label{teorema A}
Let $\X$ be a smooth hypersurface in $\Pj^n$, $n>3$, and let $\F$ be a codimension one holomorphic foliation on $\X$. If $\deg(\X)>2\deg(\F)+1$ then $\F$ extends.
\end{theorema}

We also establish an analogue criterion for the extension of codimension two distributions on smooth hypersurfaces in $\Pj^n$, $n>4$, see Theorem \ref{T:extension distribution}. The proofs of both criteria rely on the study of the restriction morphisms relating twisted differentials on $\mathbb P^n$  to twisted differentials on a projective hypersurface. Despite the simplicity of our arguments, they allow us to recover in Proposition \ref{Poincaré problem}  bounds for the degree of hypersurfaces invariant by Pfaff equations on projective spaces previously obtained in \cite{MR1452431} and \cite{MR1800822}.

In \cite{https://doi.org/10.48550/arxiv.2209.06487}, the authors showed an isomorphism between the space of foliations of degree zero on some cominuscule varieties $\X\subset \Pj^n$ and the space of foliations of degree zero on the projective space. Theorem \ref{teorema A} allows us to establish a similar result when $\X$ is a smooth projective hypersurface on $\Pj^n$, with $n>3$ and $\deg(X)>2\deg(\F)+1$.

We will show examples of non-extension in some cases not covered by Theorem \ref{teorema A}. For instance, any non-planar smooth surface $\Sp$ in $\Pj^3$ admits a foliation that does not extend. Besides, if $\Sp$ is a plane, each foliation on it has an extension, and the question of characterizing planes in $\Pj^3$ proposed by Cerveau turns out to be true.

\begin{propositione}\label{teorema B}
A smooth surface in $\Pj^3$ is a plane if, and only if, each one of its foliations extends.
\end{propositione}

We also prove the existence of a degree one foliation on the three dimensional smooth quadric, which does not extend. We will show this result using the fact that, up to a perturbation by an automorphism of $\Pj^n$, the restriction of any foliation in the projective space to a smooth hypersurface of degree at least two has an isolated singularity.

This paper is organized as follows. In Section \ref{preliminares}, we introduce Pfaff equations, holomorphic foliations, distributions, their invariant hypersurfaces, and the degree of foliations. We also define precisely the space of foliations and the restriction and extension of foliations and distributions. In Section \ref{Extension of foliation on projective hypersurfaces}, we prove Theorem \ref{teorema A} and an extension criterion for codimension two distributions. As a corollary, we obtain an isomorphism between the space of foliations on the projective space and the space of foliations on hypersurfaces, under suitable assumptions. In the last part of this section, we prove Proposition \ref{teorema B}. Finally, in Section \ref{section 5}, the property of always getting Morse singularities on restricted foliations, up to perturbation by an automorphism of the ambient projective space, is proved. This allows us to exhibit an example of foliation on a smooth quadric in $\Pj^4$ that does not extend.

\begin{akn}
{\normalfont M. G. Figueira thanks J. V. Pereira for the discussions and remarks on the results of this work. The author also thanks W. Mendson, F. Loray, G. Fazoli Domingos, G. Michels, and the anonymous referee for recommendations and comments. In particular, the close reading and detailed corrections suggested by the latter greatly improve this paper.}
\end{akn}

\section{Preliminaries}\label{preliminares}
\subsection{Foliations and Pfaff equations}  
Let $\X$ be a complex projective manifold of dimension $n$ and $L$ be a line bundle over $\X$. A {\it Pfaff equation} of codimension $q$ and coefficients in $L$ is a global section $\alpha$ of $\Omega^q_{\X}\otimes L$. The {\it singular set} of $\alpha$ is $\sing(\alpha):=\{p\in \X|\alpha(p)=0\}$.

A {\it singular codimension $q$  holomorphic foliation} $\F$ is determined by a line bundle $L$ and a Pfaff equation $\omega\in H^0(\X,\Omega^q_X\otimes L)$ such that $\codim(\sing(\omega))\geq 2$ and $\omega$ is {\it decomposable} and {\it integrable} in the following sense: let $p\in X\setminus \sing(\omega)$ be a point, then there are $1$-forms $\eta_1,\dots \eta_q$ defined over an open set $U\subset\X$ containing $p$ and  satisfying 
\begin{enumerate}[1)]
    \item $\restr{\omega}{U}=\eta_1\wedge\dots\wedge \eta_q$ (decomposability condition)
    \item $d\eta_i\wedge\eta_1\wedge\dots\wedge \eta_q=0$, for all $\,i=1,\dots,q$ (integrability condition).
\end{enumerate}
The {\it singular set of} $\F$ is $\sing(\F):=\sing(\omega)$.

\begin{rmk}
\normalfont If a divisor $D$ is contained in the singular locus of a decomposable and integrable $\omega\in H^0(\X,\Omega^q_{\X}\otimes L)$, we can replace $\omega$ with $\omega'=\frac{\omega}{f}\in H^0(\X,\Omega^q_{\X}\otimes L')$ to ensure that $D\not\subset \sing(\omega')$, where $f\in H^0(\X,\Oc_{\X}(D))$ vanishes along $D$ and $L'=L\otimes \Oc_{\X}(-D)$. This process is called {\it saturation} of $\omega$.
 \end{rmk}

For a Pfaff equation of codimension one with local representative $\eta$, the decomposability condition is automatic, and the integrability condition is given by 
$$\eta\wedge d\eta=0.$$

\begin{rmk}\label{observação sobre distribuições}
\normalfont When $\alpha\in H^0(\X,\Omega^q_{\X}\otimes L)$ satisfies only the decomposability condition, we say that $\alpha$ defines a {\it singular codimension $q$ distribution} $\D$ over $\X$. In particular, if $q=2$, then $\alpha\in H^0(\X,\Omega^q_{\X}\otimes L)$ defines distribution if and only if $\alpha\wedge \alpha=0 $ (see  \cite[Proposition 1]{https://doi.org/10.48550/arxiv.1605.09709}).
\end{rmk}

  The integrability condition ensures that the kernel of $\omega$ defines a subsheaf $T\F$ of $T\X$, called {\it tangent sheaf of} $\F$, such that in an analytic neighborhood of each non-singular point, $T\F$ is the relative tangent sheaf of a holomorphic fibration. The {\it leaves} of such a foliation are given by analytic continuation. 

\subsection{Invariant hypersurfaces}
Let $\Y\subset \X$ be a hypersurface. The inclusion map $i:\Y\rightarrow \X$ induces a restriction map of Pfaff equations projecting $\alpha\in H^0(\X,\Omega^q_{\X}\otimes L)$ on $i^*\alpha\in H^0(\Y,\Omega^q_{\Y}\otimes \restr{L}{\Y})$.

 Let $\omega\in H^0(\X, \Omega^q_{\X}\otimes L)$ be a Pfaff equation. We say that a hypersurface $\Y\subset \X$ is {\it invariant} by $\omega$ if $i^*\omega= 0$. If $\omega$ determines a foliation $\F$ and $\Y$ is invariant by $\omega$, then we say that $\Y$ is {\it invariant} by $\F$.

\subsection{Degree of a foliation}\label{degree of a foliation}
If $\Pic(\X)\simeq \Z$, taking a positive generator $M$ of $\Pic(\X)$, we define the {\it degree of a line bundle} $L$ as $\deg(L)=l$, when $L\simeq M^{\otimes l}$. If $X=\Pj^n$, $n\geq 2$, then a Pfaff equation of codimension $q$ can be represented as a degree $k$ homogeneous $q$-form of $\C^{n+1}$ such that its contraction with the radial vector field is zero, namely a $q$-form $\omega$ that satisfies $i_R(\omega)=0$, where $$R=\displaystyle\sum^{n}_{i=0}x_i\dfrac{\partial}{\partial x_i}.$$

The {\it degree of a codimension $q$ foliation} $\F$ on $\Pj^n$ is the degree of the tangency set of the leaves of $\F$ with a generic $q$-plane in $\Pj^n$, and is denoted by $\deg(\F)$.   If $\F$ is determined by $\omega\in H^0(\Pj^n,\Omega^q_{\Pj^n}(k))$, then, for instance by \cite[Lemma 3.2]{MR1800822}, one has
$$\deg(\F)=k-q-1.$$
Supposing now that $\X$ is a smooth hypersurface of $\Pj^n$, $n>3$, we have, by \cite[Corollary II.3.2]{Hartshorne1970}, $\Pic(\X)\simeq\Z$, and we define the {\it degree} of a foliation $\G$ on $\X$ generated by $\omega\in H^0(\X,\Omega^q_{\X}(k))$ as $\deg(\G)=k-q-1$.

\subsection{Restrictions and Extensions}\label{definição de extensão}
A codimension $q$ holomorphic singular foliation $\F$ (resp. a distribution $\D$) on the complex projective space $\Pj^n$ determined by $\omega\in H^0(\Pj^n,\Omega^q_{\Pj^n}(k))$ is {\it transverse} to a smooth hypersurface $\X$ if the singular set of $i^*\omega$ has codimension at least two, where $i:\X\rightarrow \Pj^n$ is the natural inclusion. In this case,  $i^*\omega$ determines a foliation $\G$ (resp. a distribution $\D'$) on $\X$, and we say that $\G$ is the {\it restriction} of $\F$ (resp. $\D'$ is the {\it restriction}). 
 
If a foliation $\G$ (resp. a distribution $\D'$) over a smooth projective hypersurface $\X\subset \Pj^n$ is the restriction of a foliation $\F$ (resp. a distribution $\D$) on $\Pj^n$ to $\X$, we say $\F$ is an {\it extension} to $\G$ (resp. an {\it extension} to $\D'$).

\subsection{Space of foliations}
 Let $X\subset \Pj^n$, $n>3$, be a smooth projective hypersurface. If $\omega\in H^0(X,\Omega^1_X(l+2))$ determines a foliation $\F$, then any non-zero constant multiple of $\omega$ determines the same foliation. Thus, we define the space of degree $l$ foliations on $\X$ as the quasi-projective variety 
     $$\Fol(\X,l):=\{[\omega]\in \Pj H^0(\X,\Omega^1_{\X}(l+2)|\,d\omega\wedge\omega=0 \text{ and }\codim(\sing(\omega))\geq 2\}.$$
 Determining the irreducible components of such a variety has been studied in some cases, especially when $\X$ is a projective space and the degree $l$ is small (see for example \cite{MR1394970} and \cite{MR4354288}).

\section{Extensions of foliations on projective hypersurfaces}\label{Extension of foliation on projective hypersurfaces}
According to Subsection \ref{definição de extensão}, in order to find an extension of a degree $l$ codimension one foliation on a projective hypersurface $\X$ determined by $\omega\in H^0(\X,\Omega^1_{\X}(l+2))$, we need to find an integrable $1$-form $\alpha\in H^0(\Pj^n,\Omega^1_{\Pj^n}(l+2))$ such that $i^*\alpha=\omega$, where $i:\X\hookrightarrow\Pj^n$ is the natural inclusion map. Thus, we need to understand the properties of the restriction map of twisted differential forms.

\begin{lemma}\label{Lema para maior ou igual a 3}
    Let $\X\subset \Pj^n$, $n>3$, be a smooth hypersurface such that $\deg(\X)\geq 2$ and $1\leq q\leq n-1$. Let
        $$\rest_q:H^0(\Pj^n,\Omega^q_{\Pj^n}(k))\rightarrow H^0(\X,\Omega^q_{\X}(k))$$
    be the restriction map of Pfaff equations. Then
        \begin{enumerate}[a)]
            \item $\rest_q$ is injective if $k-q+1\leq \deg(\X)$; 
            \item $\rest_q$ is surjective if $q<n-1$ and either $q\neq 2$ or $k\neq\deg(\X)$.
        \end{enumerate}
\end{lemma}

\begin{proof}
In \cite[Proposition 5.22]{Araujo_2018} it was proved that $\rest_1$ is an isomorphism if $k\leq \deg(\X)$ and $\rest_q$ is injective if $k-q+1\leq \deg(\X)$. Therefore, we only need to check the surjectivity of $\rest_q$ when $n-1>q$ and either $q> 2$ or $k> \deg(\X)$. For that, consider the exact sequence
\begin{equation}\label{sequência que entendo}
    0\rightarrow \Omega^q_{\Pj^n}(k-\deg(\X)) \rightarrow \Omega^q_{\Pj^n}(k)\rightarrow \restr{\Omega^q_{\Pj^n}(k)}{\X}\rightarrow 0
\end{equation} 
on $\Pj^n$ and the exact sequence
\begin{equation}\label{segunda sequência}
        0\rightarrow \Omega^{q-1}_{X}(k-\deg(\X)) \rightarrow \restr{\Omega^q_{\Pj^n}(k)}{\X}\rightarrow \Omega^q_{\X}(k)\rightarrow 0.
\end{equation}
on $\X$.

If $n-1>q>2$ or $k\neq\deg(\X)$ we have $H^1(\Pj^{n},\Omega_{\Pj^{n}}^q(k -\deg(\X)))=0$ by Bott's Theorem (see \cite[Theorem 2.3.2]{MR704986} or \cite[page 4]{okonek}). Furthermore, $H^1(\X,\Omega^{q-1}_{X}(k-\deg(\X)))=0$ by \cite[Satz 8.11]{Flenner1981}. Therefore, the maps
$$H^0(\Pj^{n},\Omega_{\Pj^{n}}^q(k))\rightarrow H^0(\X,\restr{\Omega_{\Pj^{n}}^q(k)}{\X})$$
and
$$H^0(\X,\restr{\Omega_{\Pj^{n}}^q(k)}{\X})\rightarrow H^0(\X,\Omega_{\X}^q(k))$$
obtained by long sequence in cohomology of \eqref{sequência que entendo} and \eqref{segunda sequência} are surjective. Since the restriction map is obtained by composing these two maps, the result follows. 
\end{proof}

\begin{rmk}\label{observação do caso q=2}\normalfont If $q=2$, $k=\deg(\X)$ and $\dim(\X)>3$ then $\rest_2:H^0(\Pj^n,\Omega_{\Pj^n}^2(k))\rightarrow H^0(\X,\Omega_{\X}^2(k))$ is surjective. In fact, by Bott's Theorem (\cite[Theorem 2.3.2]{MR704986}) we have the long exact sequence in cohomology of \eqref{sequência que entendo}
$$
        \begin{tikzcd}
0 \rar &  H^0(\Pj^n,\Omega^2_{\Pj^n}(k)) \rar{\phi}
             \ar[draw=none]{d}[name=X, anchor=center]{}
    & H^0(\X,\restr{\Omega^2_{\Pj^n}(k)}{\X}) \ar[rounded corners,
            to path={ -- ([xshift=2ex]\tikztostart.east)
                      |- (X.center) \tikztonodes
                      -| ([xshift=-2ex]\tikztotarget.west)
                      -- (\tikztotarget)}]{dll}[at end]{} \\      
  0 \rar & H^1(\X,\restr{\Omega^2_{\Pj^n}(k)}{\X}) \rar{\sim} & H^2(\Pj^n,\Omega^2_{\Pj^n})\simeq \C  \rar & 0.
\end{tikzcd}
$$
\normalfont Thus, $\phi$ is an isomorphism, and $H^1(\X,\restr{\Omega^2_{\Pj^n}(k)}{\X})\simeq \C$. Now \cite[5.15]{Araujo_2018} and \cite[5.17]{Araujo_2018} applied to the long sequence in cohomology of the short exact sequence (\ref{segunda sequência}) gives us 
$$        \begin{tikzcd}
0 \rar & H^0(\X,\restr{\Omega^2_{\Pj^n}(k)}{\X}) \rar{\psi}
             \ar[draw=none]{d}[name=X, anchor=center]{}
    & H^0(\X,\Omega_{\X}^2(k)) \ar[rounded corners,
            to path={ -- ([xshift=2ex]\tikztostart.east)
                      |- (X.center) \tikztonodes
                      -| ([xshift=-2ex]\tikztotarget.west)
                      -- (\tikztotarget)}]{dll}[at end]{} \\      
  H^1(\X,\Omega^1_{\X})\simeq \C \rar{\beta} & H^1(\X,\restr{\Omega_{\Pj^n}^2(k)}{X})\simeq \C \rar{} & H^1(\X,\Omega^2_{\X}(k)). &
\end{tikzcd}
$$

\normalfont If $\dim(\X)>3$ by \cite[Satz 8.11]{Flenner1981} we have $H^1(\X,\Omega_{\X}^2(k))=0$, and $\beta$ is a surjective linear application between dimension one spaces, i.e., $\beta$ is isomorphism. So $\psi$ is isomorphism, thus $\rest_2=\phi\circ \psi$ is isomorphism.
 \end{rmk}

 \begin{rmk}\label{observação sobre caso de superfícies}
     \normalfont If $n=3$ and $k\leq \deg(X)$ then $\rest_1$ is injective. Actually, if $k<\deg(\X)$ we have
     \[
        H^0(\Pj^3, \Omega^1_{\Pj^3}(k-\deg(X)))=0=H^0(\X, \Oc_{\X}(k-\deg(X )))
    \]
    by Bott's Theorem \cite[Theorem 2.3.2]{MR704986} and \cite[Theorem II.5.1 and Exercise II.5.5]{hartshorne}. So the maps $H^0(\Pj^{3},\Omega_{\Pj^{3}}^1(k))\to H^0(\X,\restr{\Omega^1_{\Pj^{3}}(k)}{\X})$ and $H^0(\X,\restr{\Omega_{\Pj^{n}}^1(k)}{\X})\to H^0(\X,\Omega_{\X}^1(k))$ are injective.

    If $\deg(\X)=k$, by Bott's Theorem \cite[Theorem 2.3.2]{MR704986} and \cite[Theorem II.5.1 and Exercise II.5.5]{hartshorne}, we have $H^1(\Pj^3,\Omega^1_{\Pj^3}(k))=0=H^1(\X,\Oc_{\X})$ and the exact sequences
   $$\begin{tikzcd}[transform canvas={scale=1}]
0 \rar &  H^0(\Pj^3,\Omega^1_{\Pj^3}(k))\rar{\phi_{1}} &  H^0(\X,\restr{\Omega^1_{\Pj^3}(k)}{\X})\rar{\alpha} & H^1(\Pj^3,\Omega^1_{\Pj^3})\simeq \C\rar &0
\end{tikzcd}$$
and
$$
\begin{tikzcd}[transform canvas={scale=1}]
0 \rar &  H^0(\X,\Oc_{\X})\simeq \C\rar{\beta} &  H^0(\X,\restr{\Omega^1_{\Pj^3}(k)}{\X})\rar{\psi_1} & H^0(\X,\Omega^1_{\X}(k))\rar &0.
\end{tikzcd}
$$

 Thus, $H^0(\Pj^3,\Omega^1_{\Pj^3}(k))$ and $H^0(\X,\Omega^1_{\X}(k)) $ are vector spaces of the same dimension. Since $\phi_1$ is injective and $\psi_1 $ is surjective, it suffices to show that the kernel of $\psi_1$ does not intersect the image of $\phi_1$ to conclude that $\rest_1:=\psi_1\circ\phi_1$ is an isomorphism.
 
 Let $f$ be a degree $k$ irreducible homogeneous polynomial such that $\X:=\{f=0\}$. We have $df\not\in H^0(\Pj^3,\Omega^1_{\Pj^3}(k))$, because $i_R(df)\neq 0$. Then, $df$  is not in the image of $\phi_1$, and $\alpha(df)$ generates $H^1(\Pj^3,\Omega^1_{\Pj^3})\simeq \C$. Furthermore, $df$ generates $\ker(\psi_1)=H^0(\X,\Oc_{\X})\simeq\C$, so $\ker(\psi_1)\cap \phi_1(H^0(\Pj^3,\Omega^1_{\Pj^3}(k)))=\emptyset$ and $\rest_1$ is an isomorphism.
 \end{rmk}

 \begin{rmk}\label{observação sobre caso de superfície 2}\normalfont If $n=3$ and $k-1\leq \deg(X)$ then $\rest_2$ is injective. In fact, by Bott's Theorem \cite[Theorem 2.3.2]{MR704986} and \cite[Satz 8.11]{Flenner1981}, we have
 \[
        H^0(\Pj^3, \Omega^2_{\Pj^3}(k-\deg(X)))=0=H^0(\X, \Omega^1_{\X}(k-\deg(X ))).
    \]
    Thus, the maps 
    $$\phi_2:H^0(\Pj^{3},\Omega_{\Pj^{3}}^2(k))\to H^0(\X,\restr{\Omega^2_{\Pj^{3}}(k)}{\X})$$
    and 
    $$\psi_2: H^0(\X,\restr{\Omega_{\Pj^{n}}^2(k)}{\X})\to H^0(\X,\Omega_{\X}^2(k))$$
    are injective, so $\rest_2:=\psi_2\circ\phi_2$ is injective.
 \end{rmk}

\begin{prop}\label{proposição transversalidade}
    Let $\X\subset\Pj^n$, $n>3$, be a smooth hypersurface, and let $\F$ be a codimension one foliation on $\Pj^n$. If $\deg(\F)+2\leq \deg(\X)$, then $\F$ is transverse to $\X$.
\end{prop}

\begin{proof}
    Let $\omega\in H^0(\Pj^n,\Omega^1_{\Pj^n}(\deg(\F)+2))$ be a $1$-form that determines $\F$. According to the foliation definition, $\codim \sing(\omega)\ge 2$.  Suppose for the sake of contradiction that there is a divisor $D\subset \sing(\rest_1(\omega))$. 
    
    As already mentioned in the general definition of degree, for $n>3$, we have $\Pic(\X)=\Z\cdot \Oc_{\X}(1)$. Therefore, $\Oc_{\X}(D) \simeq \Oc_{\X}(k)$, for some $k\in \Z_{>0}$, and there are $f\in H^0(\X,\Oc_{\X}(k))$ and $\eta\in H^0(\X,\Omega^1_{\X}(\deg(\F)+2-k))$ such that $\rest_1(\omega)=f\eta$. 

    By \cite[Exercise II.5.5]{hartshorne}, there is $g\in H^0(\Pj^n,\Oc_{\Pj^n}(k))$ whose restriction to $\X$ is $f$, and by Lemma \ref{Lema para maior ou igual a 3}, there is $\alpha\in H^0(\Pj^n,\Omega^1_{\Pj^n}(\deg(\F)+2-k))$ such that $\rest_1(\alpha)=\eta$. So,
    $$\rest_1(g\alpha)=\rest_1(\omega),$$
    and since $\deg(\F)+2\leq \deg(\X)$, Lemma \ref{Lema para maior ou igual a 3} guarantees that $\rest_1$ is injective and $g\alpha=\omega$. It follows that $\codim(\sing(\omega))=1$, contradicting $\codim \sing(\omega)\ge 2$. 
\end{proof}
 
\subsection{Proof of Theorem A}
 Let $\omega\in H^0(\X, \Omega^1_X(\deg(\F)+2))$ be a Pfaff equation determining $\F$. Let $\alpha\in H^0(\Pj^n,\Omega^1_{\Pj^n}(\deg(\F)+2))$ be an extension of $\omega$ given by Lemma \ref{Lema para maior ou igual a 3}. Now notice that the restriction of $\alpha\wedge d \alpha\in H^0(\Pj^n,\Omega^3_{\Pj^n}(2\deg(\F)+4))$ to $\X$ is zero, because $i^*\alpha=\omega$ is integrable. Again, Lemma \ref{Lema para maior ou igual a 3} guarantees the injectivity of the restriction map for 3-forms when $\deg(X)>2\deg(\F)+1$, so $\alpha\wedge d \alpha=0$ and $\F$ extends. \qed

 \subsection{Space of foliation on projective hypersurfaces} Theorem \ref{teorema A} and Proposition \ref{proposição transversalidade} allow to prove the following result.
 \begin{thm}\label{Teorema C}
    Let $\X$ be a smooth hypersurface in $\Pj^n$, $n>3$. If $\deg(\X)>2l+1$ then the map
    $$\Fol(\Pj^n,l)\rightarrow \Fol(\X,l)$$
    is an isomorphism.
\end{thm}
\begin{proof}
Let $\F$ be a codimension one foliation on $\Pj^n$ determined by $\alpha\in H^0(\Pj^n,\Omega_{\Pj^n}^1(l+2))$. Proposition \ref{proposição transversalidade} guarantees that $\F$ is transverse to $\X$, so $\rest_1(\alpha)$ defines a codimension one foliation $\G$ on $\X$ and $\deg(\G)=l$ (by very definition of the degree, see Subsection \ref{degree of a foliation}).
Therefore,  the map between quasi-projective varieties $\pi:\Fol(\Pj^n,l)\rightarrow \Fol(\X ,l)$ induced by $\rest_1$ is well-defined everywhere. It is injective, since $\deg(\X)\geq 2l+1\geq l+2$ (see Lemma \ref{Lema para maior ou igual a 3}). Furthermore, by Theorem \ref{teorema A}, it is surjective.
\end{proof}

A degree zero codimension one foliation on $\Pj^n$, $n>3$, is a pencil of hyperplanes, and it has a first integral of the form $F/G$, where $F$ and $G$ are co-prime homogeneous polynomials of degree one (for a proof of this fact, see \cite[Proposition 3.1]{MR2200857}). Thus, any two  foliations of degree zero on $\Pj^n$ are conjugated, and the set of all of them is an irreducible projective smooth variety isomorphic to the space of projective lines $\mathbb{G}(1,n)$. This fact and Theorem \ref{Teorema C} imply the following.

\begin{cor}
    Let $\X\subset\Pj^n$ be a smooth hypersurface, $n>3$. If $\deg(\X)\geq 2$, then
    $$\Fol(\X,0)\simeq \mathbb{G}(1,n).$$
\end{cor}

\subsection{Degree of invariant smooth hypersurfaces} The results regarding the injectivity of the restriction map allow us to bound the degree of smooth hypersurfaces invariant by a codimension $q$ Pfaff equation on $\Pj^n$, $n \geq 3$ and $1\leq q \leq n-1$, as follows. 

\begin{prop}\label{Poincaré problem} If a smooth hypersurface $\X\subset \Pj^n$, $n\geq 3$, is invariant by a non-trivial Pfaff equation $\alpha\in H^0(\Pj^n, \Omega^q_{\Pj^n}(k))$, $1\leq q\leq n-1$, then
$$\deg(\X)\leq k-q.$$
\end{prop}
\begin{proof}
Since $\X$ is invariant by $\alpha$, we have $\rest_q(\alpha)=0$. Suppose that $\deg(\X)\geq k-q+1$. Therefore, $\rest_q$ is injective by Lemma \ref{Lema para maior ou igual a 3}, Remark \ref{observação sobre caso de superfícies} and Remark \ref{observação sobre caso de superfície 2}. This implies $\alpha$ is zero, which is absurd.
\end{proof}

As already mentioned in the Introduction, this gives an alternative proof of previously known bounds for the degree of hypersurfaces invariant by Pfaff equations on projective spaces, see \cite{MR1452431} and \cite{MR1800822}.

\subsection{Extensions of codimension two distribution} Concerning codimension two distributions on projective smooth hypersurfaces, let us now see an example of non-extension. Let $\eta\in H^0(\Pj^4,\Omega^2_{\Pj^4}(3))$ be the codimension two Pfaff equation given by
$$ \eta=i_R(dx_0\wedge dx_1\wedge dx_2+dx_2\wedge dx_3\wedge dx_4).$$
We have $\eta\wedge \eta \neq 0$, so that $\eta$ does not satisfy the decomposability condition given in Remark \ref{observação sobre distribuições}. Therefore, it does not define a distribution on $\Pj^4$.

Let $\X\subset \Pj^4$ be a smooth hypersurface such that $\deg(\X)>2$. By Lemma \ref{Lema para maior ou igual a 3}, the restriction map 
$$\rest_2:H^0(\Pj^4,\Omega^2_{\Pj^4}(3))\rightarrow H^0(\X,\Omega^2_{\X}(3))$$
is injective. Therefore $\omega:=\rest_2(\eta)\in H^0(\X,\Omega^2_{\X}(3))$ is non-zero and satisfies $\omega \wedge \omega=0$ in $H^0(\X,\Omega^4_{\X}(6))$, since $\dim(\X)=3$. Thus $\omega$ defines distribution in $\X$. As $\rest_2^{-1}(\omega)=\{\eta\}$, the distribution determined by $\omega$ does not extend.

This non-extension example was possible because the restriction of any Pfaff equation of codimension two in $\Pj^4$ to a smooth hypersurface $\X$ automatically satisfies the decomposability condition.  On the other hand, if $\dim(\X)>3$, we have the following extension criterion.
\begin{thm}\label{T:extension distribution}
Let $\X\subset \Pj^n$, $n>4$, be a smooth hypersurface and let $\D$  be a codimension two distribution determined by $\alpha\in H^0(\X,\Omega^2_{\X}(k))$ such that $2k-3\leq \deg(\X)$. Then $\D$ extends.
\end{thm}

\begin{proof}
     If $k<3$, by \cite[Lemma 5.17]{Araujo_2018}, we have $H^0(\X,\Omega^2_{\X}(k))=0$. Suppose $k\geq 3$. Thus, by Lemma \ref{Lema para maior ou igual a 3} and Remark \ref{observação do caso q=2}, there is $\beta\in H^0(\Pj^n,\Omega^2_{\Pj^n}(k))$ whose restriction to $\X$ is $\alpha$. Since $2k-3\leq \deg(X)$, the restriction map $\rest_4: H^0(\Pj^n,\Omega^4_{\Pj^n}(2k))\rightarrow H^0(\X,\Omega^4_{\X}(2k))$ is injective. Thus $\rest_4(\beta\wedge \beta)=\alpha \wedge \alpha=0$, so $\beta\wedge \beta=0$ and, by Remark \ref{observação sobre distribuições}, $\beta$ determines a codimension two distribution on $\Pj^n$ that extends $\D$. 
\end{proof}

\subsection{Proof of Proposition \ref{teorema B}} 
Let $\Hp$ be a hyperplane in $\Pj^n$, $n\geq 3$, and $\G$ be a codimension $q$ foliation on $\Hp$. Taking a point $p\in \Pj^{n}\setminus \Hp$, we have a rational map $\pi_p:\Pj^{n}\dashrightarrow \Hp$, called  {\it projection from a point} $p$, whose definition domain is $\Pj^{n}\setminus \{p\}$. 

\begin{ex}If $\Hp=\{x_n=0\}$ and taking $p:=(0:0:\dots :0:1)\in \Pj^{n}$, $\pi_p$ is determined by
$$\begin{array}{rcl}\pi_p:\Pj^{n}&\dashrightarrow &\Hp\\
(x_0:\dots :x_{n})&\mapsto &(x_0:\dots :x_{n-1})
\end{array}.$$
\end{ex}
The pull-back of $\G$ by $\pi_p$ is a codimension $q$ foliation $\pi_p^*\G$ on $\Pj^{n}$. It is called {\it trivial extension}. Therefore, every foliation on a hyperplane of $\Pj^n$ extends. 

Let $\eta\in H^0(\Pj^3,\Omega^1_{\Pj^3}(2))$ be the contact $1$-form given by
$$
    \eta=i_R(dx_0\wedge dx_1+dx_2\wedge dx_3)=x_0dx_1-x_1dx_0+x_2dx_3-x_3dx_2.
$$
This $1$-form does not define a foliation on $\Pj^3$ because it is non-integrable. Furthermore, if $\Sp\subset \Pj^3$ is a smooth surface, with $\deg(\Sp)\geq 2$, then Proposition \ref{Poincaré problem} implies that $\Sp$ is non-invariant by $\eta$. The integrability of the restriction of $\eta$ to $\Sp$ is trivially satisfied, so that $\rest_1(\eta)$ defines a foliation on $\Sp$. 

By Remark \ref{observação sobre caso de superfícies}, if $\deg(\Sp)\geq 2$, the contact form $\eta$ is the only $1$-form whose restriction to $\Sp$ is $\rest_1(\eta)$. Therefore, the foliation generated by $\rest_1(\eta)$ on $\Sp$ does not extend to $\Pj^3$. Then, every smooth surface $\Sp$ on $\Pj^3$ such that $\deg(\Sp)\geq 2$ has a foliation that does not extend. This fact and the trivial extensions prove Proposition \ref{teorema B}.\qed

\section{A special foliation on the three-dimensional quadric}\label{section 5}

\subsection{Restriction of foliations}Let $\F$ be a codimension one foliation on $\Pj^n$ and $f\in \C[x_0,\dots ,x_n]$ be an irreducible homogeneous polynomial that determines a smooth hypersurface ${\X}$ as its  zero set. Let us give some properties on the singular set of the restriction $\restr{\F}{\X}$.

The {\it Gauss map} of $\F$ is the rational map $G_{\F}:\Pj^n\dashrightarrow \breve{\Pj}^n$ given by $p\mapsto T_p\F$, which is defined outside the singular points of $\F$.

\begin{dfn}
   Let $\G$ be a codimension one foliation on a projective hypersurface $\X\subset\Pj^{n}$. We say that $p\in\sing(\G)$ is {\it of Morse type} if there are an open subset $U$ containing $p$, and a first integral $g:U\rightarrow \C$ of $\G$ such that, in local coordinates, we can write $g=x^2_1+\dots +x^2_{n-1}$.
\end{dfn}

By Morse Lemma (see \cite[ Lemma 2.2]{book:64955}), a singular point $p$ of a function $f$ is of Morse type if and only if the {\it Hessian matrix} of $f$ in $p$, given by 
$$ \text{Hess}_p(f)=\left(\dfrac{\partial^2f(p)}{\partial x_i\partial x_j}\right)_{i,j=1,\ldots,n},$$
has a non-zero determinant. Furthermore, Morse-type singularities are isolated.

Given a generic hyperplane $\Hp\subset \Pj^n$, the singular set of $\restr{\F}{\Hp}$ is given by $(\sing(\F)\cap \Hp)\cup G_ {\F}^{-1}(\Hp)$. Furthermore, if $G_{\F}^{-1}(\Hp)$ is non-empty, then its points are isolated Morse-type singularities (see \cite[ Proposition 1.10]{ThiagoAmaralTese}).
 
 \begin{dfn} Let $X\subset \Pj^n$ be a hypersurface defined by a homogeneous polynomial $f$. The {\it Gauss map} of ${\X}$ is the rational map
  $$\begin{array}{rrcl}
  G_{\X}:&{\X}&\dasharrow &\breve{\Pj^n}\\
  &p &\mapsto &\left[\dfrac{\partial f(p)}{\partial x_0}:\dots :\dfrac{\partial f(p)}{\partial x_n}\right]
  \end{array}.$$
The indeterminacy points of $G_{\X}$ coincide with the singularities of ${\X}$.
 \end{dfn}
 
 If $\X \subset \mathbb P^n$ is a smooth projective hypersurface of degree at least two, then $G_{\X}$ is a morphism and its rank is generically equal to $n-1$ (See \cite[ Corollary I.2.5]{zak1993tangents}). The next result tells us that, up to conjugation by a generic automorphism, the restriction of a foliation $\F$ on $\Pj^n$ to a smooth hypersurface of degree at least two has singularities outside of $\sing( \F)\cap \X$ which are of Morse type.

\begin{thm}\label{Teorema das Singularidades Morse}
    Let $\F$ be codimension one holomorphic foliation on $\Pj^n$, $n>2$, and $\X\subset \Pj^n$ be a non-planar smooth hypersurface transverse to $\F$. Then there is an automorphism $h\in \Aut(\Pj^n)$ such that $\restr{h^*\F}{\X}$ has an isolated singularity $p\not\in\sing(h^*\F)$.
\end{thm}

\begin{proof} 
    Since $\X$ is smooth, the rank of the associated Gauss map is generically $n-1$. Therefore, we can suppose that the point $p=[1:0:\dots :0]\in \X$ is a regular point of $\F$, and $G_{\X}$ has maximal rank in $p$. Let us take an affine coordinate system $(x_1,\dots, x_n)$ centered on $p$ corresponding to an affine chart $U$ of $\Pj^n$ such that the first $n-1$ coordinates vectors determine the tangent space of $X$ in $p$. 

    Thus, in a neighborhood of $p$, we have that $\restr{x_n}{\X}=x_n(x_1,\dots, x_{n-1})$ vanishes in second order in $(x_1,\dots ,x_{n-1})=(0,\dots,0)$. Therefore, 
   $$x_n=\sum^{n-1}_{i,j=1}b_{ij}x_ix_j+\sum_{i\geq 3}\tilde{g_i}(x_1,\dots,x_{n-1}),$$ 
    in $X$, where each $\tilde{g_i}$ is a homogeneous polynomial of degree $i$. The rank of the matrix $(b_{ij})$ is identical to the rank of the associated Gauss map. This fact results from the description of the second fundamental form in terms of the rank of the Gauss map (see, for example, the end of Section 6.4B of \cite{https://doi.org/10.1023/A:1025366207448}). 

    Therefore, in a suitable coordinate system centered on $p$, we can write 
    $$x_n=x_1^2+\dots +x_{n-1}^2+\sum_{r\geq 3}g_i(x_1,\dots,x_{n-1}),$$ 
    in $X$, where each $g_i$ is a homogeneous polynomial of degree $i$.

    As $p$ is a regular point of $\F$, restricting $U$ if necessary, the local first integral of $\F$ defined on $U$ is of the form
    $$f=a_1x_1+\dots +a_nx_n+\sum_{j\geq 2} p_j(x_1,\dots ,x_n),$$
    where each $p_j$ is a homogeneous polynomial of degree $j$ and $a_1,\dots,a_n\in \C$, with $a_k\neq 0$ for some $k\in\{1,\dots, n\}$.  Then for every $\lambda\in \C^*$, we can find an automorphism $h_{\lambda}\in \text{Aut}(\Pj^n)$ such that the local first integral of $h_{\lambda}^*\F$ on an open set $V\subset U$ of $p$ is
    $$f_{\lambda}=\lambda \cdot x_n+\sum_{j\geq 2} h_j(x_1,\dots ,x_n),$$
    where $h_j$ is a homogeneous polynomial of degree $j$ and $h_2=\displaystyle\sum^{n}_{j=1}\displaystyle\sum^{n}_{i=1}a_{ij}x_ix_j$.
    For instance, take $h_{\lambda}$ the automorphism that preserves the hyperplane at infinity and, on affine coordinates, is defined by $x_i\mapsto x_i$, if $i\neq k$, and 
    $$x_k\mapsto\frac{1}{a_k}\left(\left(\sum_{\substack{i=1 \\ i\neq k}}^{n}-a_ix_i\right)+\lambda x_n\right).$$
     
     The restriction of $h_{\lambda}^*\F$ to $\X$ has a first integral on $p$ of the form
     $$\restr{f_{\lambda}}{\X}= \lambda(x_1^2+\dots +x^2_{n-1}) + \sum_{j\geq 2} \tilde{ h_j}(x_1,\dots ,x_{n-1}),$$
     such that each $\tilde{h_j}$ is a homogeneous polynomial of degree $j$. Additionally, 
     $$\tilde{h_2}=h_2(x_1,\dots, x_{n-1},0)=\displaystyle\sum^{n-1}_{j=1}\displaystyle\sum^{n- 1}_{i=1}a_{ij}x_ix_j,$$
     and the $a_{ij}=a_{ji}\in \C$ do not depend on the choice of $\lambda$. In fact, if $\lambda_1,\lambda_2\in \C^*$, to replace $f_{\lambda_1}$ by $f_{\lambda_2}$ just make a pull-back of $h^*_{\lambda_1}\F$ via an automorphism of $\Pj^n$ that perform  $x_i\mapsto x_i$, if $i\neq n$, and $x_n\mapsto \frac{\lambda_2 x_n} {\lambda_1}$, in affine coordinates.
    
    As $\dfrac{\partial \restr{f_{\lambda}}{\X}}{\partial x_i}(p)=0$ for all $i=1,\dots ,n-1$, then $p$ is a singularity of $\restr{h_{\lambda}^*\F}{\X}$. We will show that for generic $\lambda\in \C^*$, $p$ is a Morse singularity, so it is an isolated singularity. To do this, it suffices to show that, for a generic $\lambda$, the determinant of the Hessian matrix of $\restr{f_{\lambda}}{\X}$ is nonzero at the point $p$.
    
    In fact, we have
    $$\dfrac{\partial^2\restr{f_{\lambda}}{\X}}{\partial x_i^2}(p)=2\lambda+2a_{ii}$$
    and, if $i\neq j$, then
    $$\dfrac{\partial^2 \restr{f_{\lambda}}{\X}}{\partial x_i\partial x_j}(p)=2a_{ij}.$$
    
    Therefore, the determinant of the Hessian matrix $$\text{Hess}_p(\restr{f_{\lambda}}{\X})=\det\left(\dfrac{\partial^2\restr{f_{\lambda}}{\X}}{\partial x_i\partial x_j}(p)\right)$$
    is a polynomial of degree $n-1$ in the variable $\lambda$, with leading coefficient equal to $2^{n-1}$. Then $\text{Hess}_p(\restr{f_{\lambda}}{\X})$ is not identically zero and has finite number of roots. This tells us that for a generic $\lambda\in \C^*$, $p$ is a Morse-type singularity of $\restr{h_{\lambda}^*\F}{\X}$.
\end{proof}

\subsection{Non-extension of a foliation on the three-dimensional quadric}
For a codimension one and degree one foliation over a smooth quadric in $\Pj^4$, Theorem \ref{teorema A} does not guarantee its extension. In \cite[5.11]{https://doi.org/10.48550/arxiv.1112.3871}, the authors presented a degree one foliation over the smooth quadric on $\Pj^4$ with interesting properties. We will briefly describe such foliation and prove that this foliation does not extend.

We can realize $\Pj^4$ as the equivalence class of homogeneous polynomials of degree four in two variables, i.e., four unordered points on the projective line. With this identification, there is a natural action of $\Aut(\Pj^1)\simeq\PGL(2,\C)$ on $\Pj^4$. The closure of the orbit of $\PGL(2,\C)$ on the point of $\Pj^4$ corresponding to $\{1,-1,i,-i\}$ generates a smooth quadric $\Qu^ 3\subset \Pj^4$.

The action of the affine subgroup $\mathbb{A}\text{ff}(\C)\subset\PGL(2,\C)$ on $\Qu^3$ induces a codimension one foliation $\mathcal{A}$ of degree one, whose singular set does not contain isolated points: $\sing(\mathcal{A})$ is composed of a rational normal curve of degree four, a Veronese curve of degree 3 and a line, which correspond, respectively, to points of the form $\{p,p,p,p\}$, $\{\infty,p,p,p\}$ and $\{p,\infty,\infty,\infty\}$. Also, $\mathcal{A}$ belongs to an irreducible component $\mathbb{A}\text{ff}\subset\Fol(\Qu^3,1)$ called {\it affine component} whose general element is given by conjugation of $\mathcal{A}$ with element of $ \Aut(\Qu^3)\simeq \PO(5,\C)$ (see \cite[Theorem 5.2]{https://doi.org/10.48550/arxiv.1112.3871}).

\begin{thm}
   Let $\mathcal{A}$ be the degree one affine foliation described above. Then $\mathcal{A}$ does not extend.
\end{thm}

\begin{proof}
    Suppose by contradiction that there is a foliation $\F$ in $\Pj^4$ which is the extension of $\mathcal{A}$. We have the following rational map
    $$\begin{array}{rcl}
       \Phi:\Aut(\Pj^4) &\dashrightarrow &\Fol(Q,1)  \\
         g              &\mapsto         &\restr{g^*\F}{Q}  
    \end{array}.$$
    As $\Phi(id)=\mathcal{A}$ and $\Aut(\Pj^4)$ is irreducible, then component $\mathbb{A}\text{ff}$ contains the image of $\Phi$.
    
    On the other hand, by Theorem \ref{Teorema das Singularidades Morse}, there exists an automorphism $h\in \Aut(\Pj^4)$ such that $\restr{h^*\F} {\Qu^3}$ has an isolated singularity $p\not\in\sing(h^*\F)$. Therefore, such foliation can not belong to $\mathbb{A}\text{ff}$, as the elements in this component are a pull-back of $\mathcal{A}$ by an automorphism of $Q$ and does not have isolated singularities, which is absurd.
\end{proof}

\bibliography{references}{}
\bibliographystyle{amsplain}
\end{document}